\documentclass[a4paper,twoside,12pt]{article}
\usepackage{amssymb} 
\begin{document}
\title{Lower bound on the number of the maximum genus embedding of $K_{n,n}$ \footnotemark[2]
\author{Guanghua Dong$^{1,2}$, Han Ren$^{3}$, Ning Wang$^{4}$, Yuanqiu Huang$^{1}$\\
{\small\em 1.Department of Mathematics, Normal
University of Hunan, Changsha, 410081, China}\\
\hspace{-1mm}{\small\em 2.Department of Mathematics, Tianjin
Polytechnic
University, Tianjin, 300387, China}\\
\hspace{-1mm}{\small\em 3.Department of Mathematics, East China
Normal University, Shanghai, 200062,China}\\
\hspace{-5mm}{\small\em 4.Department of Information Science and
Technology,
Tianjin University of Finance }\\
\hspace{-74mm} {\small\em and Economics, Tianjin, 300222, China}\\
}} \footnotetext[2]{\footnotesize \em   This work was partially
Supported  by the China Postdoctoral Science Foundation funded
project (Grant No: 20110491248), the National Natural Science
Foundation of China (Grant No: 11171114), and  the New Century
Excellent Talents in University (Grant No: NCET-07-0276).}
 \footnotetext[1]{\footnotesize \em E-mail: gh.dong@163.com(G. Dong). }

\date{}
\maketitle

\vspace{-.8cm}
\begin{abstract}
In this paper, we  provide an method to obtain the lower bound on
the number of the distinct maximum genus embedding of the complete
bipartite graph $K_{n,n}$ ($n$ be an odd number), which, in some
sense, improves the results of S. Stahl and H. Ren.

\end{abstract}


\bigskip
\noindent {\bf 1. Introduction}

\bigskip
Graphs considered here are all connected and finite. A $surface$ $S$
means a compact and connected two-manifold without boundaries. A
$cellular$ $embedding$ of a graph $G$ into a surface $S$ is a
one-to-one mapping $\psi:$ $G \rightarrow S$ such that each
component of $S-\psi(G)$ is homomorphic to an open disc. The maximum
genus $\gamma_M(G)$ of a connected graph \emph{G} is the maximum
integer \emph{k} such that there exists an embedding of $G$ into the
orientable surface of genus $k$. By Euler's polyhedron formula, if a
cellular embedding of a graph $G$ with $n$ vertices, $m$ edges and
$r$ faces is on an orientable surface of genus $\gamma$, the
$n-m+r=2-2\gamma$. Since $\gamma \geqslant 1$, we have $\gamma(G)
\leqslant \frac{1}{2}\lfloor\beta(G)\rfloor$, where $\beta(G)=m-n+1$
is called the $Betti$ $number$ (or $cycle$ $rank$) of the graph $G$.
It follows that $\gamma_M(G)\leqslant
\frac{1}{2}\lfloor\beta(G)\rfloor$. If $\gamma_M(G)=
\frac{1}{2}\lfloor\beta(G)\rfloor$, then the graph is called $upper$
$embeddable$. It is not difficult to deduced that a graph is upper
embeddable if and only if its face number is not greater than two.
Since the introductory investigations on the maximum genus of graphs
by Nordhaus, Stewart, and White$^{\cite{nor}}$, this parameter has
attracted considerable attention from mathematicians and computer
scientists. Up to now, the research about the maximum genus of
graphs mainly focus on the aspects as characterizations and
complexity, the upper embeddability, the lower bound, the
enumeration of the distinct maximum genus embedding, $etc.$. For
more detailed information, the reader can be found in a survey  in
$\cite{top}$.

It is well known that the enumeration of the distinct maximum genus
embedding plays an important role in the study of the genus
distribution problem, which  may be used to decide whether two given
graphs are isomorphic. It was S. Stahl$^{\cite{sta}}$ who provides
the first result about the lower bound on the number of the distinct
maximum genus embedding, which is states as the following:

\medskip

{\bf Lemma 1$^{\cite{sta}}$} \  \ A connected graph (loops and
multi-edges are allowed) of order $n$ with degree sequence $d_1$,
$d_2$, $\dots$, $d_{n}$ has at least
$$(d_1-5)!(d_2-5)!(d_3-5)!(d_4-5)!\prod_{i=5}^{n}(d_{i}-2)!$$
distinct orientable embeddings with at most two facial walks, where
$m!=1$ whenever $m\leqslant0$.

\medskip

But up to now, except \cite{sta} and \cite{ren}, there is little
result concerning the number of the maximum genus embedding of
graphs. In this paper, we will provide a method to enumerate the
number of the distinct maximum genus embedding of the complete
bipartite graph $K_{n,n}$ ($n$ be an odd number), and offer a lower
bound which is better than that of S. Stahl$^{\cite{sta}}$ and H.
Ren$^{\cite{ren}}$ in some sence. Furthermore, the enumerative
method below can be used to any maximum genus embedding, other than
the method in \cite{sta} which is restricted to upper embeddable
graphs. Terminologies and notations not explained here can be seen
in $\cite{bon}$ for general graph theory, and in $\cite{liu}$ and
$\cite{moh}$ for topological graph theory.

\bigskip

\noindent {\bf 2. Main results}

\bigskip

A simple graph $G$ is called a $complete$ $bipartite$ $graph$ if its
vertex set can be partitioned into two subsets $X$ and $Y$ so that
every edge has one end in $X$ and one end in $Y$, and every vertex
in $X$ is joined to every vertex in $Y$. We denote a $complete$
$bipartite$ $graph$ $G$ with bipartition $X$ and $Y$  by
$G_{[X][Y]}$. A 2-$path$ is called a $v$-$type$-$edge$, and is
denoted by $\mathcal {V}$. Let $\psi(G)$ be an embedding of a graph
$G$. We say that a $v$-$type$-$edge$ are inserted into $\psi(G)$ if
the three endpoints of the $v$-$type$-$edge$ are inserted into the
corners of the faces in $\psi(G)$, yielding an embedding of
$G+\mathcal {V}$. The embedding $\psi(G)$ of $G$ is called a
$one$-$face$-$embedding$ (or $two$-$face$-$embedding$) if the total
face number of $\psi(G)$ is one (or two). The following observation
can be easily obtained and is essential in the proof of the Theorem
A.

\medskip

{\bf Observation} \  \ Let $\psi(G)$ be an embedding of a graph $G$.
We can insert a $v$-$type$-$edge$ $\mathcal {V}$ to $\psi(G)$ to get
an embedding  $\rho(G+\mathcal {V})$ of $G+\mathcal {V}$ so that the
face number of $\rho(G+\mathcal {V})$ is not more than that of
$\psi(G)$.

\medskip

{\bf Theorem A} \  \ For $n\equiv1 \ (mod \ 2)$, the number of the
distinct maximum genus embedding of the complete bipartite graph
$K_{n,n}$ is at least
\begin{displaymath}
2^{\frac{n-1}{2}}\times \big((n-2)!! \big)^{n}\times\big( (n-1)!
\big)^{n}.
\end{displaymath}

\medskip

{\bf Proof} \ \  Let $n=2s+1$ and $V(K_{n,n})=\{x_1, x_2, \dots,
x_{n}\}\cup \{y_1, y_2, \dots, y_{n}\}$, where $X=\{x_1, x_2, \dots,
x_{n}\}$ and $Y=\{y_1, y_2, \dots, y_{n}\}$ are the two independent
set of $K_{n,n}$. We denote the $v$-$type$-$edge$
$y_{2i}x_{j}y_{2i+1}$ by $\mathcal {V}_{ji}$,  where $i\in \{1,2,
\dots, s\}$ and $j\in\{1,2, \dots, n\}$.

\bigskip

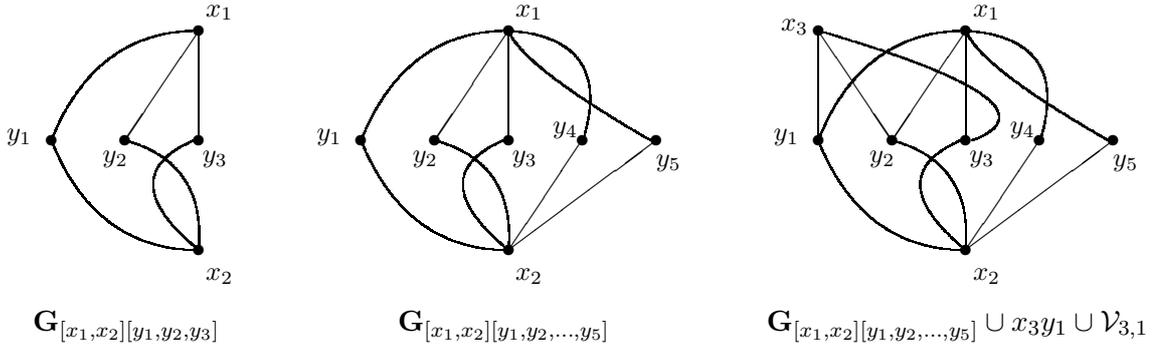
\setlength{\unitlength}{0.97mm}
\begin{center}
\begin{picture}(100,40)

\put(-20,5) {\begin{picture}(10,10)

\put(0,20){\circle*{1.5}}   

\put(10,20){\circle*{1.5}}   

\put(20,20){\circle*{1.5}}   

\put(20,35){\circle*{1.5}}   

\put(20,5){\circle*{1.5}}   

\qbezier(20,35)(6,35)(0,20)  

\put(20,35){\line(-2,-3){10}}   

\put(20,35){\line(0,-1){15}}  

\qbezier(20,5)(6,5)(0,20)  

\qbezier(20,5)(21,17)(10,20)   

\qbezier(20,5)(8,15)(19.3,20)  

\begin{footnotesize}

\put(-6,20){{$y_1$}}

\put(7,17){{$y_2$}}

\put(20.5,17){{$y_3$}}

\put(21,37){{$x_1$}}

\put(21,1){{$x_2$}}

\end{footnotesize}

\begin{small}

\put(-2.5,-6){{\bf G$_{[x_1, x_2][y_1, y_2, y_3]}$}}

\end{small}

\end{picture}}


\put(22,5) {\begin{picture}(10,10)

\put(0,20){\circle*{1.5}}   

\put(10,20){\circle*{1.5}}   

\put(20,20){\circle*{1.5}}   

\put(30,20){\circle*{1.5}}   

\put(40,20){\circle*{1.5}}   

\put(20,35){\circle*{1.5}}   

\put(20,5){\circle*{1.5}}   

\qbezier(20,35)(6,35)(0,20)  

\qbezier(20,35)(34,34)(30.3,20.7)  

\put(20,35){\line(-2,-3){10}}   

\qbezier(20,35)(22,30)(39.3,20)  

\put(20,35){\line(0,-1){15}}  

\qbezier(20,5)(8,15)(19.3,20)  

\qbezier(20,5)(6,5)(0,20)  

\put(20,5){\line(2,3){9.6}}   

\qbezier(20,5)(21,17)(10,20)   

\put(20,5){\line(4,3){19.5}}  

\begin{footnotesize}

\put(-6,20){{$y_1$}}

\put(7,17){{$y_2$}}

\put(20.5,17){{$y_3$}}

\put(26,21){{$y_4$}}

\put(40,16.5){{$y_5$}}

\put(21,37){{$x_1$}}

\put(21,1){{$x_2$}}

\end{footnotesize}

\begin{small}

\put(5,-6){{\bf G$_{[x_1, x_2][y_1, y_2, \dots, y_5]}$}}

\end{small}

\end{picture}}


\put(84,5){\begin{picture}(10,10)

\put(0,35){\circle*{1.5}}   

\put(0,35){\line(0,-1){15}}   

\put(0,20){\circle*{1.5}}   

\put(0,35){\line(2,-3){10}}  

\put(10,20){\circle*{1.5}}   

\qbezier(0,35)(34,25)(20.9,20)  

\put(20,20){\circle*{1.5}}   

\put(30,20){\circle*{1.5}}   

\put(40,20){\circle*{1.5}}   

\put(20,35){\circle*{1.5}}   

\put(20,5){\circle*{1.5}}   

\qbezier(20,35)(6,35)(0,20)  

\qbezier(20,35)(34,34)(30.3,20.7)  

\put(20,35){\line(-2,-3){10}}   

\qbezier(20,35)(22,30)(39.3,20)  

\put(20,35){\line(0,-1){15}}  

\qbezier(20,5)(8,15)(19.3,20)  

\qbezier(20,5)(6,5)(0,20)  

\put(20,5){\line(2,3){9.6}}   

\qbezier(20,5)(21,17)(10,20)   

\put(20,5){\line(4,3){19.5}}  

\begin{footnotesize}

\put(-6,20){{$y_1$}}

\put(7,17){{$y_2$}}

\put(20.5,17){{$y_3$}}

\put(26,21){{$y_4$}}

\put(40,16.5){{$y_5$}}

\put(21,37){{$x_1$}}

\put(21,1){{$x_2$}}

\put(-5,35){{$x_3$}}

\end{footnotesize}

\begin{small}

\put(-7,-6){{\bf G$_{[x_1, x_2][y_1, y_2, \dots, y_5]}\cup
x_3y_1\cup\mathcal {V}_{3,1}$}}

\end{small}

\end{picture}}

\end{picture}
\end{center}

\medskip

{\bf Claim 1:}  \  \ For $G_{[x_1, x_2][y_1, y_2, \dots, y_{n}]}$,
the number of the distinct $one$-$face$-$embedding$ is at least
$2^{s}\times((2s-1)!!)^2$.

\medskip

There are 2 different ways to embed $G_{[x_1, x_2][y_1, y_2,
y_{3}]}$ on an orientable surface so that the embedding is a
$one$-$face$-$embedding$. Select any one of them and denote its face
boundary by $W_0$. In $W_0$, there are three $face$-$corner$
containing $x_1$ and $x_2$ respectively. So, there are 3 different
ways to put $\mathcal {V}_{1,2}$ in $W_0$, and 3 different ways to
put $\mathcal {V}_{2,2}$ in $W_0$. Therefore, the total number of
ways to put $\mathcal {V}_{1,2}\cup\mathcal {V}_{2,2}$ in $W_0$ is
$3\times3=9$. For each of the above 9 ways, there are 2 different
ways to make the embedding of $G_{[x_1, x_2][y_1, y_2, \dots,
y_{5}]}$ being a $one$-$face$-$embedding$. So, for each of the
$one$-$face$-$embedding$ of $G_{[x_1, x_2][y_1, y_2, y_{3}]}$, there
are $3\times3\times2$ different ways to add $\mathcal {V}_{1,2}\cup
\mathcal {V}_{2,2}$ to $G_{[x_1, x_2][y_1, y_2, y_{3}]}$ to get a
$one$-$face$-$embedding$ of $G_{[x_1, x_2][y_1, y_2, \dots,
y_{5}]}$.

Similarly, we can get that for each of the $one$-$face$-$embedding$
of $G_{[x_1, x_2][y_1, y_2, \dots, y_{5}]}$, there are
$5\times5\times2$ different ways to add $\mathcal {V}_{1,3}\cup
\mathcal {V}_{2,3}$ to $G_{[x_1, x_2][y_1, y_2, \dots, y_{5}]}$ to
get a $one$-$face$-$embedding$ of $G_{[x_1, x_2][y_1, y_2, \dots,
y_{7}]}$.

In general, we have that for each of the $one$-$face$-$embedding$ of
$G_{[x_1, x_2][y_1, y_2, \dots, y_{2k-1}]}$, there are
$(2k-1)\times(2k-1)\times2$ different ways to add $\mathcal
{V}_{1,k}\cup \mathcal {V}_{2,k}$ to $G_{[x_1, x_2][y_1, y_2, \dots,
y_{2k-1}]}$ to get a $one$-$face$-$embedding$ of $G_{[x_1, x_2][y_1,
y_2, \dots, y_{2k+1}]}$.

From the above we can get that the number of the distinct
$one$-$face$-$embedding$ of $G_{[x_1, x_2][y_1, y_2, \dots, y_{n}]}$
is at least
\begin{eqnarray*}
2\times(3\times3\times2)\times(5\times5\times2)\times(7\times7\times2)\times
\dots \times((2s-1)\times(2s-1)\times2)  \\
\lefteqn{ = 2^{s}\times((2s-1)!!)^2. }  \hspace*{141mm} \\
\end{eqnarray*}

\vspace{-7mm}

 {\bf Claim 2:}  \  \ For each of the
$one$-$face$-$embedding$ of $G_{[x_1, x_2][y_1, y_2, \dots,
y_{n}]}$, there are at least $2\times(2s-1)!!\times2^{2s}$ different
ways to make $G_{[x_1, x_2, x_3][y_1, y_2, \dots, y_{n}]}$ being  a
$one$-$face$-$embedding$.

\medskip

Let $\mathcal {E}_1$ be an arbitrary $one$-$face$-$embedding$ of
$G_{[x_1, x_2][y_1, y_2, \dots, y_{n}]}$.  In $\mathcal {E}_1$,
there are two different $face$-$corner$ containing $y_{i} \
(i=1,2,3)$. So, there are $2\times2\times2(=8)$ different ways to
add $y_1x_3\cup\mathcal {V}_{3,1}$ to $\mathcal {E}_1$ to make
$G_{[x_1, x_2][y_1, y_2, \dots, y_{n}]}\cup y_1x_3\cup\mathcal
{V}_{3,1}$ being a $one$-$face$-$embedding$. For each of the above 8
$one$-$face$-$embedding$ of $G_{[x_1, x_2][y_1, y_2, \dots,
y_{n}]}\cup y_1x_3\cup\mathcal {V}_{3,1}$, there are 3 different
$face$-$corner$ containing $x_3$ and 2 different $face$-$corner$
containing $y_{i} \ (i=4,5)$. So, for each of the above 8
  $one$-$face$-$embedding$ of $G_{[x_1, x_2][y_1, y_2, \dots,
y_{n}]}\cup y_1x_3\cup\mathcal {V}_{3,1}$, there are
$3\times2\times2$ different ways to add $\mathcal {V}_{3,2}$ to
$G_{[x_1, x_2][y_1, y_2, \dots, y_{n}]}\cup y_1x_3\cup\mathcal
{V}_{3,1}$  to  make $G_{[x_1, x_2][y_1, y_2, \dots, y_{n}]}\cup
y_1x_3\cup\mathcal {V}_{3,1}\cup\mathcal {V}_{3,2}$ being a
$one$-$face$-$embedding$.

In general, we have that for each of the $one$-$face$-$embedding$ of
$G_{[x_1, x_2][y_1, y_2, \dots, y_{n}]}\cup y_1x_3\cup\mathcal
{V}_{3,1}\cup\mathcal {V}_{3,2}\cup\dots\cup\mathcal {V}_{3,k-1}$,
there are $(2k-1)\times2\times2$ different ways to add $\mathcal
{V}_{3,k}$ to $G_{[x_1, x_2][y_1, y_2, \dots, y_{n}]}\cup
y_1x_3\cup\mathcal {V}_{3,1}\cup\mathcal
{V}_{3,2}\cup\dots\cup\mathcal {V}_{3,k-1}$ to get a
$one$-$face$-$embedding$ of $G_{[x_1, x_2][y_1, y_2, \dots,
y_{n}]}\cup y_1x_3\cup\mathcal {V}_{3,1}\cup\mathcal
{V}_{3,2}\cup\dots\cup\mathcal {V}_{3,k-1}\cup\mathcal {V}_{3,k}$.

From the above we can get that for each of the
$one$-$face$-$embedding$ of $G_{[x_1, x_2][y_1, y_2, \dots,
y_{n}]}$, there are at least
\begin{eqnarray*}
(2\times2\times2)\times(3\times2\times2)\times(5\times2\times2)\times
\dots \times((2s-1)\times2\times2)  \\
\lefteqn{ = 2\times(2s-1)!!\times2^{2s} }  \hspace*{121mm} \\
\end{eqnarray*}

\vspace{-7mm} \hspace{-8mm} different ways to make $G_{[x_1, x_2,
x_3][y_1, y_2, \dots, y_{n}]}$ being  a $one$-$face$-$embedding$.
\medskip

{\bf Claim 3:}  \  \ For each of the $one$-$face$-$embedding$ of
$G_{[x_1, x_2, x_3][y_1, y_2, \dots, y_{n}]}$, there are at least
$3\times(2s-1)!!\times3^{2s}$ different ways to make $G_{[x_1, x_2,
x_3, x_4][y_1, y_2, \dots, y_{n}]}$ being  a
$one$-$face$-$embedding$.

\medskip

Let $\mathcal {E}_2$ be an arbitrary $one$-$face$-$embedding$ of
$G_{[x_1, x_2, x_3][y_1, y_2, \dots, y_{n}]}$.  In $\mathcal {E}_2$,
there are three different $face$-$corner$ containing $y_{i} \
(i=1,2,3)$. So, there are $3\times3\times3(=27)$ different ways to
add $y_1x_4\cup\mathcal {V}_{4,1}$ to $\mathcal {E}_2$ to make
$G_{[x_1, x_2, x_3][y_1, y_2, \dots, y_{n}]}\cup y_1x_4\cup\mathcal
{V}_{4,1}$ being a $one$-$face$-$embedding$. For each of the above
27 $one$-$face$-$embedding$ of $G_{[x_1, x_2, x_3][y_1, y_2, \dots,
y_{n}]}\cup y_1x_4\cup\mathcal {V}_{4,1}$, there are 3 different
$face$-$corner$ containing $x_4$ and 3 different $face$-$corner$
containing $y_{i} \ (i=4,5)$. So, for each of the above 27
  $one$-$face$-$embedding$ of $G_{[x_1, x_2, x_3][y_1, y_2, \dots,
y_{n}]}\cup y_1x_4\cup\mathcal {V}_{4,1}$, there are
$3\times3\times3$ different ways to add $\mathcal {V}_{4,2}$ to
$G_{[x_1, x_2,x_3][y_1, y_2, \dots, y_{n}]}\cup y_1x_4\cup\mathcal
{V}_{4,1}$  to  make $G_{[x_1, x_2, x_3][y_1, y_2, \dots,
y_{n}]}\cup y_1x_4\cup\mathcal {V}_{4,1}\cup\mathcal {V}_{4,2}$
being a $one$-$face$-$embedding$.

In general, we have that for each of the $one$-$face$-$embedding$ of
$G_{[x_1, x_2, x_3][y_1, y_2, \dots, y_{n}]}\cup y_1x_4\cup\mathcal
{V}_{4,1}\cup\mathcal {V}_{4,2}\cup\dots\cup\mathcal {V}_{4,k-1}$,
there are $(2k-1)\times3\times3$ different ways to add $\mathcal
{V}_{4,k}$ to $G_{[x_1, x_2, x_3][y_1, y_2, \dots, y_{n}]}\cup
y_1x_4\cup\mathcal {V}_{4,1}\cup\mathcal
{V}_{4,2}\cup\dots\cup\mathcal {V}_{4,k-1}$ to get a
$one$-$face$-$embedding$ of $G_{[x_1, x_2, x_3][y_1, y_2, \dots,
y_{n}]}\cup y_1x_4\cup\mathcal {V}_{4,1}\cup\mathcal
{V}_{4,2}\cup\dots\cup\mathcal {V}_{4,k-1}\cup\mathcal {V}_{4,k}$.

From the above we can get that for each of the
$one$-$face$-$embedding$ of $G_{[x_1, x_2, x_3][y_1, y_2, \dots,
y_{n}]}$, there are at least
\begin{eqnarray*}
(3\times3\times3)\times(3\times3\times3)\times(5\times3\times3)\times
\dots \times((2s-1)\times3\times3)  \\
\lefteqn{ = 3\times(2s-1)!!\times3^{2s} }  \hspace*{121mm} \\
\end{eqnarray*}

\vspace{-7mm} \hspace{-8mm} different ways to make $G_{[x_1, x_2,
x_3, x_4][y_1, y_2, \dots, y_{n}]}$ being  a
$one$-$face$-$embedding$.
\medskip

Similarly, we can get the following general result.

\medskip

{\bf Claim 4:}  \  \ For each of the $one$-$face$-$embedding$ of
$G_{[x_1, x_2, \dots,  x_{k-1}][y_1, y_2, \dots, y_{n}]}$, there are
at least $(k-1)\times(2s-1)!!\times(k-1)^{2s}$ different ways to
make $G_{[x_1, x_2, \dots,  x_{k-1},  x_{k}][y_1, y_2, \dots,
y_{n}]}$ being  a $one$-$face$-$embedding$.

\medskip

Noticing that a $one$-$face$-$embedding$ of a graph must be its
maximum genus embedding, we can get, from Claim 1 - Claim 4, that
the number of the distinct maximum genus embedding of $K_{n,n}$ is
at least
\begin{eqnarray*}
\{2^{s}\times((2s-1)!!)^2\}\times\{2\times(2s-1)!!\times2^{2s}\}\times\{3\times(2s-1)!! \\
\lefteqn{ \times3^{2s}\}\times\dots \times\{2s\times(2s-1)!!\times(2s)^{2s}\} }  \hspace*{109mm} \\
\lefteqn{ = 2^{s}\times((2s-1)!!)^{2s+1}\times((2s)!)^{2s+1}}  \hspace*{114mm} \\
\lefteqn{ = 2^{\frac{n-1}{2}}\times((n-2)!!)^{n}\times((n-1)!)^{n}.    \hspace*{75mm} \Box}  \hspace*{114mm} \\
\end{eqnarray*}

\vspace{-4mm}

{\bf Remark}   \ \ Through a comparison we can get that  the result
in Theorem A is much better than that of Lemma 1$^{\cite{sta}}$ when
$n\leqslant9$.

In \cite{ren}, the second author of the present paper obtained that
a connected loopless graph of order $n$ has at least
$\frac{1}{4^{\gamma_{M}(G)}}\prod_{v\in V(G)}(d(v)-1)!$ distinct
maximum genus embedding. Let $f_{1}(n)= 2^{\frac{n-1}{2}}\times\big(
(n-2)!! \big)^{n}\times \big( (n-1)! \big)^{n}$,  $f_2(n)=
\frac{1}{4^{\gamma_{M}(G)}}\prod_{v\in V(G)}\big(d(v)-1
\big)!=\frac{1}{4^{\frac{(n-1)(n-1)}{2}}}\times\big((n-1)!
\big)^{2n} $. Through a computation we can get $f_1(3) - f_2(3)=16$,
$f_1(5) - f_2(5)=6772211712$. So, when $n\leqslant5$ the result
obtained in Theorem A is much better than that of \cite{ren}.


\end{document}